\documentclass[12pt]{article}

 \usepackage{a4}                          
 \usepackage{amsmath}
 \usepackage{amssymb}
 \newtheorem{thm}{Theorem}[section]

 \newtheorem{cor}[thm]{Corollary}

 \newtheorem{lemma}[thm]{Lemma}
 \newtheorem{rem}[thm]{Remark}
 \DeclareMathOperator{\Log}{Log}
 \DeclareMathOperator{\Arg}{Arg}
\renewcommand{\t}{\theta}
 \newcommand{\eps}{\varepsilon}

 \title{A Pick function  related to the sequence of volumes of the
   unit ball in $n$-space\footnote{Both authors acknowledge support by grant
272-07-0321 from the Danish Research Council for Nature and Universe.}}

 \author{Christian Berg $^{\dagger}$ and Henrik L. Pedersen
   $^{\ddagger}$\\
 \footnotesize $\dagger$ \ Institute of Mathematical Sciences,
   University of Copenhagen\\
    \footnotesize Universitetsparken 5; DK-2100 K\o benhavn \O, Denmark\\
\footnotesize E-mail berg@math.ku.dk 
\\  \footnotesize
   $\ddagger$ \ Department of Basic Sciences and Environment\\
\footnotesize Faculty of Life Sciences, University of Copenhagen\\
\footnotesize Thorvaldsensvej 40, DK-1871 Frederiksberg C\\
\footnotesize E-mail henrikp@dina.kvl.dk
}

 \date{\today}

\begin{document}

 \maketitle

 \begin{abstract}
 We show that 
 $$
 F_a(x)=\frac{\ln \Gamma (x+1)}{x\ln(ax)}
 $$
is a Pick function for $a\ge 1$ and find its integral representation.
 We also consider the function
$$
f(x)=\left(\frac{\pi^{x/2}}{\Gamma(1+x/2)}\right)^{1/(x\ln x)}
$$
and show that $\ln f(x+1)$ is a Stieltjes function and that $f(x+1)$ is
completely monotonic on $]0,\infty[$. In particular
$f(n)=\Omega_n^{1/(n\ln n)},n\ge 2$ is a
  Hausdorff moment sequence. Here $\Omega_n$ is the volume of the unit
  ball in Euclidean $n$-space.
 \end{abstract}

\noindent 
2010 {\em Mathematics Subject Classification}:
primary 33B15; secondary 30E20, 30E15.

\noindent
Keywords: gamma function, completely monotone function.

 \section{Introduction and results}
 Since the appearance of the paper \cite{A:V:V}, monotonicity
 properties
 of the functions
\begin{equation}\label{eq:Fa}
 F_a(x)=\frac{\ln \Gamma (x+1)}{x\ln(ax)}, \quad x>0, a>0
 \end{equation}
 have attracted the attention of several authors in connection with
 monotonicity properties of the volume $\Omega_n$ of the unit ball in Euclidean
 $n$-space. A recent paper about inequalities
 involving $\Omega_n$ is \cite{Al}. 

Let us first consider the case $a=1$. In \cite{B:P1} the authors proved
that $F_1$ is a {\it Bernstein function}, which means that it is positive and
has a completely monotonic derivative, i.e.,   
\begin{equation}\label{eq:BP1}
 (-1)^{n-1}F_1^{(n)}(x)\geq 0,\quad x>0, n\ge 1.
\end{equation}
 This extended monotonicity and concavity proved in \cite{A:Q} and
 \cite{E:L} respectively.
 
 We actually proved a stronger statement than \eqref{eq:BP1}, namely that the reciprocal function
 $x\ln x/\ln\Gamma(x+1)$ is a Stieltjes transform, i.e. belongs to the
  Stieltjes cone $\mathcal S$ of functions of the form
  \begin{equation}\label{eq:int-rep}
  g(x)=c+\int_0^\infty\frac{d\mu(t)}{x+t},\quad x>0,
  \end{equation}
  where $c\geq 0$ and $\mu$ is a non-negative measure on $[0,\infty[$ satisfying
  $$\int_0^\infty\frac{d\mu(t)}{1+t}<\infty.$$

The result was obtained using the holomorphic extension of the
function $F_1$ to the cut plane $\mathcal A=\mathbb C\setminus]-\infty,0]$,
leading to an explicit formula for the measure $\mu$ in
\eqref{eq:int-rep}. Our derivation used the fact that the holomorphic function $\log\Gamma(z)$
only vanishes in $\mathcal A$ at the points $z=1$ and $z=2$, a result
interesting in itself and included as an appendix in \cite{B:P1}. A
simpler proof of the non-vanishing of $\log\Gamma(z)$ appeared in \cite{B:P2}.

In a subsequent paper \cite{B:P2} we proved an almost equivalent 
result, namely that $F_1$ is a Pick function, and obtained the following
representation formula
\begin{equation}\label{eq:F1}
F_1(z)=1- \int_0^\infty \frac{d_1(t)}{z+t}\,dt,\quad z\in\mathcal A
\end{equation}
where 
\begin{equation}\label{eq:d}
d_1(t)=\frac{\ln|\Gamma(1-t)|+(k-1)\ln t}{t((\ln
  t)^2+\pi^2)}\quad\mbox{for}
\quad t\in \left]k-1,k\right[,\quad k=1,2,\ldots
\end{equation}
and $d_1(t)$ tends to infinity when $t$ approaches $1,2,\ldots$. Since
$d_1(t)>0$ for $t>0$, \eqref{eq:BP1} is an immediate consequence of \eqref{eq:F1}.

 We recall that a Pick function is holomorphic function $\varphi$ in the upper
half-plane $\mathbb H=\{z=x+iy \in\mathbb C\mid y>0\}$ satisfying
$\Im\varphi(z)\ge 0$ for $z\in\mathbb H$, cf. \cite{D}.

For $a=2$ Anderson and Qiu proved in \cite{A:Q} that $F_2$ is
strictly increasing on $[1,\infty[$, thereby proving a conjecture from
\cite{A:V:V}. Alzer proved in \cite{Al} that $F_2$ is concave on 
$[46,\infty[$. In \cite{Q:G} the concavity was extended to the optimal
interval $]\tfrac12,\infty[$.

We will now describe the main results of the present paper.

We also denote by $F_a$ the holomorphic extension of \eqref{eq:Fa} to
$\mathcal A$ with an isolated singularity at $z=1/a$, which is
 a simple pole with residue $\ln\Gamma(1+1/a)$ assuming $a\neq 1$, while $z=1$ is a
removable singularity for $F_1$. For details about this extension see
the beginning of section 2. Using the residue theorem we obtain:

\begin{thm}\label{thm:1} For $a>0$ the function $F_a$ has the integral representation 
\begin{equation}\label{eq:Farep}
F_a(z)=1+\frac{\ln\Gamma(1+1/a)}{z-1/a}-\int_0^\infty\frac{d_a(t)}{z+t}\,dt,
\quad z\in\mathcal A\setminus\{1/a\},
\end{equation}
where 
\begin{equation}\label{eq:darep}
d_a(t)=\frac{\ln|\Gamma(1-t)|+(k-1)\ln(at)}{t((\ln
  (at))^2+\pi^2)}\quad \mbox{for}
\quad t\in \left]k-1,k\right[,\quad k=1,2,\ldots,
\end{equation}
and $d_a(0)=0, d_a(k)=\infty,k=1,2,\ldots$.
We have $d_a(t)\ge 0$ for $t\ge 0,a\ge 1/2$\footnote{This is slightly
  improved in Remark~\ref{thm:bettera} below.} and $F_a$ is a Pick function
for $a\ge 1$ but not for $0<a<1$.
\end{thm}

From this follows the monotonicity property conjectured in \cite{Q:G}:
\begin{cor}\label{thm:cor1} Assume $a\ge 1$. Then 
\begin{equation}\label{eq:Faprop}
(-1)^{n-1}F_a^{(n)}(x)>0,\quad x>1/a,n=1,2,\ldots.
\end{equation}
\end{cor}
In particular, $F_a$ is strictly increasing and strictly concave on
the interval $]1/a,\infty[$.

The function
\begin{eqnarray}\label{eq:vol}
f(x)=\left(\frac{\pi^{x/2}}{\Gamma(1+x/2)}\right)^{1/(x\ln x)}
\end{eqnarray}
has been studied because the volume $\Omega_n$ of the unit ball  in
 $\mathbb R^n$ is
$$
\Omega_n=\frac{\pi^{n/2}}{\Gamma(1+n/2)},n=1,2,\ldots.
$$

We prove the following integral representation of the extension of
$\ln f(x+1)$ to the cut plane $\mathcal A$. 

\begin{thm}\label{thm:2} For $z\in\mathcal A$ we have
\begin{equation}\label{eq:logf}
\log
f(z+1)=-\frac12+\frac{\ln(2/\sqrt{\pi})}{z}+\frac{\ln(\sqrt{\pi})}{\Log(z+1)}+\frac12\int_1^\infty\frac{d_2((t-1)/2)}{z+t}\,dt.
\end{equation}
In particular $1/2+\log f(x+1)$ is a Stieltjes function and $f(x+1)$
is completely monotonic.
\end{thm}

We recall that completely monotonic functions
$\varphi:\left]0,\infty\right[\to\mathbb R$ are characterized by
Bernstein's theorem as
\begin{equation}\label{eq:Bern}
\varphi(x)=\int_0^\infty e^{-xt}\,d\mu(t),
\end{equation}
where $\mu$ is a positive measure on $[0,\infty[$ such that the
integrals above make sense for all $x>0$.

We also recall that a sequence $\{a_n\}_{n\ge 0}$ of positive numbers is a
Hausdorff moment sequence if it has the form
\begin{equation}\label{eq:Hau}
a_n=\int_0^1 x^n\,d\sigma(x),\;n\ge 0,
\end{equation} 
where $\sigma$ is a positive measure on the unit interval. Note that $\lim_{n\to\infty}a_n=\sigma(\{1\})$.
For a discussion of these concepts see \cite{B:C:R} or \cite{W}. It is
clear that if $\varphi$ is completely monotonic with the integral
representation \eqref{eq:Bern}, then
$a_n=\varphi(n+1),n\ge 0$ is a Hausdorff moment sequence, because
$$
a_n=\int_0^\infty e^{-(n+1)t}\,d\mu(t)=\int_0^1 x^n\,d\sigma(x),
$$
where $\sigma$ is the image measure of $e^{-t}\,d\mu(t)$ under
$e^{-t}$.
Since $\lim_{x\to\infty}f(x+1)=e^{-1/2}$ we get

\begin{cor}\label{thm:cor2}
 The sequence
\begin{equation}\label{eq:f(n)}
f(n+2)=\Omega_{n+2}^{1/((n+2)\ln(n+2))}, n=0,1,\ldots
\end{equation}
is a Hausdorff moment sequence tending to $e^{-1/2}$.
\end{cor}

A Hausdorff moment sequence is clearly decreasing and convex and by
the Cauchy-Schwarz inequality is is even logarithmically convex,
meaning that $a_n^2\le a_{n-1}a_{n+1},\;n\ge 1$. The
latter property
was obtained in \cite{Q:G} in a different way.

\section{Properties of the function $F_a$}

In this section we will study the holomorphic extension of the function $F_a$ defined in
\eqref{eq:Fa}. First a few words about notation. We use $\ln$ for the
natural logarithm but only applied to positive numbers. The
holomorphic extension of $\ln$ from the open half-line $]0,\infty[$ to
the cut plane $\mathcal A=\mathbb C\setminus ]-\infty,0]$ is denoted
$\Log z=\ln|z|+i\Arg z$, where $-\pi<\Arg z<\pi$ is the principal
argument. The holomorphic branch of the logarithm of $\Gamma(z)$ for $z$ in the simply
connected domain $\mathcal A$ which equals $\ln\Gamma(x)$ for $x>0$ is
denoted $\log\Gamma(z)$. The imaginary part of $\log\Gamma(z)$ is a
continuous branch of argument of $\Gamma(z)$ which we denote
$\arg\Gamma(z)$, i.e.,
$$
\log\Gamma(z)=\ln|\Gamma(z)|+i\arg\Gamma(z),\;z\in\mathcal A.
$$
We shall use the following property of
$\log\Gamma(z)$, cf. \cite[Lemma 2.1]{B:P1}

 \begin{lemma}\label{thm:logGamma}
 We have, for any $k\geq 1$,
 $$
 \lim_{z\to t,\Im z>0}\log \Gamma (z)= \ln|\Gamma (t)| -i\pi k
 $$
 for $t\in ]-k,-k+1[$ and
 $$
 \lim_{z\to t,\Im z>0}|\log \Gamma (z)|=\infty
 $$
 for $t=0,-1,-2, \ldots$.
 \end{lemma}

The expression
$$
F_a(z)=\frac{\log\Gamma(z+1)}{z\Log(az)}
$$
clearly defines a holomorphic function in $\mathcal A\setminus\{1/a\}$,
and $z=1/a$ is a simple pole unless $a=1$, where the residue
$\ln\Gamma(1+1/a)$ vanishes.

\begin{lemma}\label{thm:behonR} 
For $a>0$ and $t \le 0$ we have 
\begin{equation}\label{eq:ytozero}
\lim_{y\to 0^+} \Im F_a(t+iy)=\pi d_a(-t),
\end{equation}
where $d_a$ is given by \eqref{eq:darep}.
\end{lemma}

{\it Proof.} For $-1<t<0$ we get 
$$
\lim_{y\to 0^+} F_a(t+iy)=\frac{\ln\Gamma(1+t)}{t(\ln(a|t|)+i\pi)},
$$
hence $\lim_{y\to 0^+}\Im F_a(t+iy)=\pi d_a(-t)$.
For $-k<t<-k+1,\, k=2,3,\ldots$ we find using Lemma~\ref{thm:logGamma}
 $$
\lim_{y\to 0^+}
F_a(t+iy)=\frac{\ln|\Gamma(1+t)|-i(k-1)\pi}{t(\ln(a|t|)+i\pi)},
$$
hence $\lim_{y\to 0^+}\Im F_a(t+iy)=\pi d_a(-t)$ also in this case. 

For $t=-k,\;k=1,2,\ldots$ we have 
$$
|F_a(-k+iy)|\ge
\frac{\left|\ln|\Gamma(-k+1+iy)|\right|}{|-k+iy||\Log(a(-k+iy))|}\to\infty
$$
for $y\to 0^+$ because $\Gamma(z)$ has poles at
$z=0,-1,\ldots$. Finally, for $t=0$ we get \eqref{eq:ytozero} from the
next Lemma.
$\quad\square$

\begin{lemma}\label{thm:behatzero} For $a>0$ we have
$$
\lim_{z\to 0,z\in\mathcal A}|F_a(z)|=0.
$$
\end{lemma}

{\it Proof.} Since $\log\Gamma(z+1)/z$ has a removable singularity for
$z=0$ the result follows because
$|\Log(az)|\ge |\ln(a|z|)|\to \infty$ for $|z|\to 0,z\in\mathcal A$.$\quad\square$

\begin{lemma}\label{thm:behoncircles} For $a>0$ we have the radial
  behaviour
\begin{equation}\label{eq:radial}
\lim_{r\to\infty} F_a(re^{i\t})=1\;\mbox{for}\; -\pi<\t<\pi,
\end{equation}
and there exists a  constant $C_a>0$
 such that for $k=1,2,\ldots$ and $-\pi<\t<\pi$
\begin{equation}\label{eq:bhc1}
|F_a((k+\tfrac12)e^{i\t})|\le C_a.
\end{equation}
\end{lemma}

{\it Proof.} We first note that
\begin{equation}\label{eq:quot}
F_a(z)=F_1(z)\frac{\Log(z)}{\Log(az)},
\end{equation}
and since
$$
\lim_{|z|\to\infty,z\in\mathcal A}\frac{\Log(z)}{\Log(az)}=1
$$
it is enough to prove the results for $a=1$. We do this by using a
method introduced in \cite[Prop. 2.4]{B:P1}.

 Define 
$$
R_k=\{ z = x+iy \in \mathbb C \,\mid \,-k\leq x < -k+1,\,
0<y\leq 1\, \}\;\;\mbox{for}\;\; k\in \mathbb Z
$$
 and
 $$
R=\cup_{k=0}^{\infty}R_k,\quad S=\{ z = x+iy \in \mathbb C \,\mid
\,x\le 1,|y|\le 1\}.
$$

By Lemma~\ref{thm:logGamma} it is clear that 
\begin{equation}\label{eq:Mk}
M_k=\sup_{|\t|<\pi}|F_1((k+\tfrac12)e^{i\t})|<\infty
\end{equation}
for each $k=1,2,\ldots$, so it is enough to prove that $M_k$ is
bounded for $k\to\infty$.

 Stieltjes (\cite[formula 20]{S}) found the following
formula for $\log \Gamma (z)$ for $z$ in the cut plane $\mathcal{A}$
\begin{equation}\label{eq:St1}
\log \Gamma (z+1)= \ln\sqrt{2\pi} + (z+1/2)\Log z -z +\mu(z).
\end{equation}
Here
$$
\mu(z)=\sum_{n=0}^{\infty} h(z+n)= \int_0^{\infty} \frac{P(t)}{z+t}dt,
$$
where $h(z)=(z+1/2)\Log (1+1/z) -1$ and $P$ is periodic with period
1 and $P(t)=1/2-t$ for $t\in [0,1[$. A derivation of these formulas can
also be found in \cite{Ar}. The integral above is improper,
and integration by parts yields
\begin{equation}\label{eq:mu}
\mu(z)=\frac{1}{2} \int_0^{\infty} \frac{Q(t)}{(z+t)^2}dt,
\end{equation}
where $Q$ is periodic with period 1 and $Q(t)=t-t^2$ for $t\in [0,1[$.
Note that by \eqref{eq:mu} $\mu$ is a completely monotonic
function. For further properties of Binet's function $\mu$ see \cite{K:P}.

We claim that
$$
|\mu(z)|\leq \frac{\pi}{8} \;\mbox{for}\; z \in \mathcal
A\setminus S.
$$
 In fact, since $0\leq Q(t)\leq 1/4$, we get for $z=x+iy
\in \mathcal A$
$$
|\mu(z)|\leq \frac{1}{8}\int_0^{\infty}\frac{dt}{(t+x)^2+y^2}.
$$
For $x>1$ we have
$$
\int_0^{\infty}\frac{dt}{(t+x)^2+y^2}
\leq \int_0^{\infty}\frac{dt}{(t+1)^2} =1,
$$
and for $x\le 1, |y|\geq 1$ we have
$$
\int_0^{\infty}\frac{dt}{(t+x)^2+y^2}
= \int_x^{\infty}\frac{dt}{t^2+y^2}  < \int_{-\infty}^\infty \frac{dt}{t^2+1}=\pi.
$$
Since
$$
F_1(z) =
1+ \frac{ \ln\sqrt{2\pi} + 1/2\Log z -z +\mu(z)}{z\Log z},
$$
for $z \in \mathcal A$, we immediately get \eqref{eq:radial} and 
\begin{equation}
\label{eq:upper1}
|F_1(z)|\le 2
\end{equation}
for all  $z\in \mathcal A\setminus S$ for which $|z|$ is sufficiently
large. In particular, there exists  $N_0\in\mathbb N$ such that
\begin{equation}\label{eq:bhc2}
|F_1((k+\tfrac12)e^{i\t})|\le 2\; 
\mbox{for}\; k\ge N_0,\;(k+\tfrac12)e^{i\t}\in\mathcal A\setminus S.
\end{equation}

By continuity the quantity
\begin{equation}\label{eq:cont}
c=\sup\left\{ |\log\Gamma(z)|\;\mid\;z=x+iy,\tfrac12\le x\le 1,0\le y \le 1\right\}
\end{equation}
is finite.

We will now estimate the quantity $|F_1((k+\tfrac12)e^{i\t})|$ when
$(k+\tfrac12)e^{i\t}\in S$, and since
$F_1(\overline{z})=\overline{F_1(z)}$, it is enough to consider the
case when $(k+\tfrac12)e^{i\t}\in R_{k+1}$. To do this we use the
relation
 \begin{equation}
 \label{eq:smart}
 \log \Gamma (z+1) = \log \Gamma (z+k+1)-\sum_{l=1}^{k} \Log (z+l)
 \end{equation}
 for $z \in \mathcal A$ and $k\in\mathbb N$. Equation \eqref{eq:smart} follows from the
 fact that the functions on both sides of the equality sign are holomorphic
 functions in $\mathcal A$, and they agree on the positive half-line by
 repeated applications of the functional equation for the Gamma function.

For $z=(k+\tfrac12)e^{i\t}\in R_{k+1}$ we get $|\log\Gamma(z+k+1)|\le
c$ by \eqref{eq:cont}, and hence by \eqref{eq:smart}
$$
|\log\Gamma(z+1)|\le c+\sum_{l=1}^k |\Log(z+l)|\le
c+k\pi+\sum_{l=1}^k |\ln|z+l||.
$$
For $l=1,\ldots,k-1$ we have $k-l<|z+l|<k+2-l$, hence
$0<\ln|z+l|<\ln(k+2-l)$. Furthermore, $1/2\le |z+k|\le \sqrt{2}$,
hence $-\ln 2<\ln|z+k|\le (\ln 2)/2$. Inserting this we get
$$
|\log\Gamma(z+1)|\le c+k\pi +\sum_{j=2}^{k+1} \ln j <
c+k\pi+k\ln(k+1).
$$
From this we get for $z=(k+\tfrac12)e^{i\t}\in R_{k+1}$ 
\begin{equation}\label{eq:upper2}
|F_1(z)|\le \frac{c+k\pi+k\ln(k+1)}{(k+\tfrac12)\ln(k+\tfrac12)}
\end{equation}
which tends to 1 for $k\to\infty$. Combined with \eqref{eq:bhc2} we
see that there exists $N_1\in\mathbb N$ such that
$$
|F_1((k+\tfrac12)e^{i\t})|\le 2\;\mbox{for}\; k\ge N_1,\,-\pi<\t<\pi,
$$
which shows that $M_k$ from \eqref{eq:Mk} is a bounded sequence.
$\quad\square$

\begin{lemma}\label{thm:boundFa} Let $a>0$. For $k=1,2,\ldots$ there exists an
 integrable function
 $f_{k,a}:\left]-k,-k+1\right[\to\left[0,\infty\right]$ such that
\begin{equation}\label{eq:boundFa}
|F_a(x+iy)|\le f_{k,a}(x)\;\mbox{for}\; -k<x<-k+1,0<y\le 1.
\end{equation}
\end{lemma}

{\it Proof.} For $z=x+iy$ as above we get using \eqref{eq:smart}
$$
|\log\Gamma(z+1)|\le |\log\Gamma(z+k+1)|+\sum_{l=1}^k|\Log(z+l)|
\le L+k\pi +\sum_{l=1}^k|\ln|z+l||,
$$
where $L$ is the maximum of $|\log\Gamma(z)|$ for
$z\in\overline{R_{-1}}$.
We only treat the case $k\ge 2$ because the case $k=1$ is a simple
modification  combined with Lemma~\ref{thm:behatzero}.

For $l=1,\ldots,k-2$ we have $1<|z+l|<1+k-l$,  and for $l=k-1,k$
$\ln|x+l|\le \ln|z+l|\le (1/2)\ln 2$, so we find
\begin{equation}\label{eq:boundFa1}
|\log\Gamma(z+1)|\le L+k\pi+\sum_{j=2}^k \ln j+|\ln|x+k-1||+|\ln|x+k||,
\end{equation}
so as $f_{k,1}$ we can use the right-hand side of \eqref{eq:boundFa1}
divided by $(k-1)\ln(k-1)$. Using \eqref{eq:quot} we next define
$$
f_{k,a}(x)=f_{k,1}(x)\max_{z\in\overline{R_k}}\frac{|\Log
  z|}{|\Log(az)|}.
$$
  $\quad\square$ 

\medskip
{\bf Proof of Theorem~\ref{thm:1}}
For fixed $w\in\mathcal A\setminus\{1/a\}$ we choose $\eps>0,k\in\mathbb N$ such that
$\eps < |w|,1/a < k+\tfrac12$ and consider the positively oriented contour
$\gamma(k,\eps)$ in
$\mathcal A$ consisting of the half-circle $z=\eps e^{i\t},\t\in
[-\tfrac{\pi}2,\tfrac{\pi}2]$ and the half-lines $z=x \pm i\eps,x\le 0$
until they cut the  circle $|z|=k+\tfrac12$, which closes the contour.
By the residue theorem we find
$$
\frac{1}{2\pi i}\int_{\gamma(k,\eps)}
\frac{F_a(z)}{z-w}\,dz=F_a(w)+\frac{\ln\Gamma(1+1/a)}{1/a-w}.
$$

We now let $\eps\to 0$ in the contour integration. By
Lemma~\ref{thm:behatzero} the contribution from the half-circle with
radius $\eps$ will tend to zero, and by Lemma~\ref{thm:behonR} and Lemma~\ref{thm:boundFa} we get
$$
\frac{1}{2\pi}\int_{-\pi}^{\pi}\frac{F_a((k+\tfrac12)e^{i\t})}{(k+\tfrac12)e^{i\t}-w}(k+\tfrac12)e^{i\t}\,d\t
+\int_{-k-\tfrac12}^0 \frac{d_a(-t)}{t-w}\,dt=F_a(w)+\frac{\ln\Gamma(1+1/a)}{1/a-w}.
 $$

For $k\to\infty$ the integrand in the first integral converges to 1
for each $\t\in\left]-\pi,\pi\right[$ and by
Lemma~\ref{thm:behoncircles} Lebesgue's theorem on dominated
convergence can be applied, so we finally get
$$
F_a(w)=1+\frac{\ln\Gamma(1+1/a)}{w-1/a}-\int_0^\infty\frac{d_a(t)}{t+w}\,dt.
$$

The last integral above appears as an improper integral, but we shall
see that the integrand is
Lebesgue integrable. We show below that $d_a(t)\ge 0$ when $a\ge 1/2$
and for these values of $a$ the integrability is obvious. The function
$d_a$ tends to 0 for $t\to 0$ and  has a logarithmic singularity at
$t=1$ so $d_a$ is integrable over $]0,1[$. For $k-1<t<k,\;k\ge 2$ we
have
\begin{equation}\label{eq:intda}
d_a(t)=\frac{(\ln(t))^2+\pi^2}{(\ln(at))^2+\pi^2}d_1(t)+\frac{(k-1)\ln
  a}{t\left((\ln(at))^2+\pi^2\right)},
\end{equation}
and the factor in front of $d_1(t)$ is a bounded continuous function
with limit 1 at 0 and at infinity. Therefore 
$$
\int_1^\infty \frac{|d_a(t)|}{t}\,dt<\infty
$$
follows from the finiteness of the corresponding integral for $a=1$
provided that we establish
$$
S:=\sum_{k=2}^\infty
(k-1)\int_{k-1}^k\frac{dt}{t^2\left((\ln(at))^2+\pi^2\right)}<\infty.
$$
Choosing $N\in\mathbb N$ such that $aN>1$, we can estimate

\begin{eqnarray*}
S &<& \sum_{k=1}^\infty \int_{ka}^{(k+1)a}\frac{dt}{t(\ln^2(t)+\pi^2)}
<\int_a^{Na}\frac{dt}{t(\ln^2(t)+\pi^2)}+\sum_{k=N}^\infty\int_{ka}^{(k+1)a}\frac{dt}{t\ln^2(t)}\\
&=&\int_a^{Na}\frac{dt}{t(\ln^2(t)+\pi^2)}+\frac{1}{\ln(aN)}<\infty.
\end{eqnarray*}

We next examine positivity of $d_a$.

For $0<t<1$ we have
$$
d_a(t)=\frac{\ln|\Gamma(1-t)|}{t((\ln(at))^2+\pi^2)}>0
$$
because $\Gamma(s)>1$ for $0<s<1$.

For $k\ge 2$ and $t\in\left]k-1,k\right[$ the numerator $N_a$ in $d_a$ can
be written
$$
N_a(t)=\ln\Gamma(k-t)+\sum_{l=1}^{k-1}\ln\frac{ta}{t-l},
$$
where we have used the functional equation for $\Gamma$, hence
$$
N_a(t) \ge \sum_{l=1}^{k-1}\ln\frac{k}{k-l}+(k-1)\ln a=(k-1)\ln
k-\ln\Gamma(k)+(k-1)\ln a,
$$
because $\Gamma(k-t)>1$ and $t/(t-l)$ is decreasing for
$k-1<t<k$. From \eqref{eq:St1} we get
\begin{equation}\label{eq:lnGamma}
\ln\Gamma(k)=\ln\sqrt{2\pi} +(k-1/2)\ln k -k +\mu(k)
\end{equation}
and in particular for $k=2$
$$
\mu(2)=2-\frac{3}{2}\ln 2-\ln\sqrt{2\pi}.
$$
Using \eqref{eq:lnGamma} we find
$$
N_a(t)\ge k - \frac12\ln k -\ln\sqrt{2\pi}-\mu(k) +(k-1)\ln a \ge  k -
\frac12\ln k -2+\frac{3}{2}\ln 2 +(k-1)\ln a,
$$
because $\mu$ is decreasing on $]0,\infty[$ as shown by \eqref{eq:mu}.

For $a\ge 1/2$ and $k-1<t<k$ with $k\ge 2$ we then get
$$
N_a(t)\ge k(1-\ln 2)-\frac12\ln k +\frac{5}{2}\ln 2-2\ge 0,
$$
because the sequence $c_k,k\ge 2$ on the right-hand side is increasing with $c_2=0$.

 We also see that $d_a(t)$ tends to infinity for $t$ approaching the
 end points of the interval $]k-1,k[$. For $z=1/a+iy,y>0$ we get from
 \eqref{eq:Farep}
$$
\Im F_a(1/a+iy)=-\frac{\ln\Gamma(1+1/a)}{y}+\int_0^\infty\frac{yd_a(t)}{(1/a+t)^2+y^2}\,dt.
$$
The last term tends to 0 for $y\to 0$ while the first term tends to
$-\infty$ when $0<a<1$. This shows that $F_a$ is not a Pick function
for these values of $a$.
$\quad\square$

\begin{rem}\label{thm:bettera} {\rm We proved in Theorem~\ref{thm:1}
    that $d_a(t)$ is non-negative on $[0,\infty[$ for $a\ge 1/2$. This
    is not best possible, and we shall explain that the smallest value of $a$ for which
    $d_a(t)$ is non-negative is $a_0=0.3681154742..$.

Replacing $k$ by $k+1$ in the numerator $N_a$ for $d_a$ given by
\eqref{eq:darep}, we see that
$$
N_a(t)=\ln|\Gamma(1-t)|+k\ln(at) \;\;\mbox{for}\;\; t\in
]k,k+1[,\;\;k=1,2,\ldots
$$
is non-negative if and only if
$$
\ln(1/a) \le \ln(k+s)+\frac{1}{k}\ln|\Gamma(1-k-s)| \;\;\mbox{for}\;\; s\in
]0,1[,\;\;k=1,2,\ldots,
$$
and using the reflection formula for $\Gamma$ this is equivalent to
$\ln(1/a)\le \rho(k,s)$ for all $0<s<1$ and all $k=1,2,\ldots$, where
\begin{equation}\label{eq:rho}
\rho(k,s)= \ln (k+s)-\frac{1}{k}\ln\left(\Gamma(k+s)\frac{\sin(\pi
    s)}{\pi}\right).
\end{equation}
Using Stieltjes' formula \eqref{eq:St1}, we find that 
\begin{eqnarray}\label{eq:necs}
\lefteqn{\rho(k,s)= 1+\frac{\ln(\pi/2)}{2k}}\nonumber\\
&&-(1/k)\left[(s-1/2)\ln(s+k)+\ln\sin(\pi s)
-s+\mu(s+k)\right]
\end{eqnarray}
for all $s\in\left]0,1\right[$ and $k=1,2,\ldots$.
For fixed $s\in \left]0,1\right[$ we see that $\rho(k,s)\to 1$ as
 $k\to\infty$, so
$\ln(1/a) \le 1$ is a necessary condition for non-negativity of
$d_a(t)$. This condition is not sufficient, because for $\ln(1/a)=1$
the inequality $1\leq \rho(k,s)$
 is equivalent to
$$
0\ge (1/2)\ln(2/\pi)+(s-1/2)\ln(s+k)+\ln\sin(\pi s)
-s+\mu(s+k)
$$
which does not hold when $k$ is sufficiently large and  $1/2<s<1$. 

For each $k=1,2,\ldots$ it is easy to verify that the function $\rho_k(s)=\rho(k,s)$ has a unique minimum
$m_k$ over $]0,1[$, and clearly 
\begin{equation}\label{eq:smallest}
\ln(1/a_0)=\inf\{m_k,k\ge 1\}
\end{equation}
 determines the
smallest value of $a$ for which $d_a(t)$ is non-negative. Using Maple
one obtains that $m_k$ is decreasing for $k=1,\ldots,510$ and increasing
for $k\ge 510$ with limit 1. Therefore $m_{510}=\inf m_k=0.9993586013..$ corresponding to $a_0=0.3681154742..$. 
We add that $m_1=1.6477352344..,m_{178}=1.0000028637.. ,m_{179}=0.9999936630..$.}
\end{rem}

\section{Properties of the function $f$}

{\bf Proof of Theorem~\ref{thm:2}} The function
$$
\ln f(x)=\frac{(x/2)\ln \pi-\ln\Gamma(1+x/2)}{x\ln x}
$$
clearly has a meromorphic extension to $\mathcal A\setminus{1}$ with a
simple pole at $z=1$ with residue $\ln 2$. We denote this meromorphic
extension $\log f(z)$ and have
$$
\log f(z+1)=\frac{\ln\sqrt{\pi}}{\Log(z+1)}-\frac12 F_2\left(\frac{z+1}2\right).
$$
Using the representation \eqref{eq:Farep}, we immediately get
\eqref{eq:logf}. It is well-known that $1/\Log(z+1)$ is a Stieltjes
function, cf. \cite[p.130]{B:F}, and the integral representation is
\begin{equation}\label{eq:log}
\frac{1}{\Log(z+1)}=\int_1^\infty\frac{dt}{(z+t)((\ln(t-1))^2+\pi^2)}.
\end{equation}
It follows that $\ln(\sqrt{e}f(x+1))$ is a Stieltjes function, in
particular completely monotonic, showing that $\sqrt{e}f(x+1)$ belongs
to the class $\mathcal L$ of logarithmically completely monotonic
functions studied in \cite{F:G:C} and in \cite{B1}. Therefore also
$f(x+1)$ is completely monotonic.$\quad\square$ 

\section{Representation of $1/F_a$}

For $a>0$ we consider the function 
\begin{equation}\label{eq:Ga}
G_a(z)=1/F_a(z)=\frac{z\Log(az)}{\log\Gamma(z+1)}
\end{equation}
which is holomorphic in $\mathcal A$ with an isolated singularity at
$z=1$, which is a simple pole with residue $\ln a/\Psi(2)=\ln
a/(1-\gamma)$ if $a\ne 1$, while it is a removable singularity when
$a=1$. Here $\Psi(z)=\Gamma'(z)/\Gamma(z)$ and $\gamma$ is Euler's constant.  

\begin{thm}\label{thm:Ga} For $a>0$ the function $G_a$ has the integral representation 
\begin{equation}\label{eq:Garep}
G_a(z)=1+\frac{\ln a}{(1-\gamma)(z-1)}+\int_0^\infty\frac{\rho_a(t)}{z+t}\,dt,
\quad z\in\mathcal A\setminus\{1\},
\end{equation}
where 
\begin{equation}\label{eq:rhoarep}
\rho_a(t)=t\frac{\ln|\Gamma(1-t)|+(k-1)\ln(at)}{(\ln
  |\Gamma(1-t)|)^2+((k-1)\pi)^2}\quad \mbox{for}
\quad t\in \left]k-1,k\right[,\quad k=1,2,\ldots,
\end{equation}
and $\rho_a(0)=1/\gamma, \rho_a(k)=0,\;k=1,2,\ldots$, which makes
$\rho_a$ continuous on $[0,\infty[$.
We have $\rho_a(t)\ge 0$ for $t\ge 0,\;a\ge a_0=0.3681154742..$, cf. Remark~\ref{thm:bettera}, and $G_a(x+1)$ is a Stieltjes function
for $a\ge 1$ but not for $0<a<1$.
\end{thm}

{\it Proof.} We notice that
for $-k<t<-k+1,\, k=1,2,\ldots$ we get using Lemma~\ref{thm:logGamma}
 $$
\lim_{y\to 0^+}
G_a(t+iy)=\frac{t(\ln(a|t|)+i\pi)}{\ln|\Gamma(1+t)|-i(k-1)\pi},
$$
and for $t=-k,k=1,2,\ldots$ we get
 $$
\lim_{y\to 0^+}
|G_a(-k+iy)|=0
$$
because of the poles of $\Gamma$,
hence $\lim_{y\to 0^+}\Im G_a(t+iy)=-\pi \rho_a(-t)$ for $t<0$. 

For fixed $w\in\mathcal A\setminus\{1\}$ we choose $\eps>0,k\in\mathbb N$ such that
$\eps < |w|,1 < k+\tfrac12$ and consider the positively oriented contour
$\gamma(k,\eps)$ in $\mathcal A$ which was used in the proof of Theorem~\ref{thm:1}.

By the residue theorem we find
$$
\frac{1}{2\pi i}\int_{\gamma(k,\eps)}
\frac{G_a(z)}{z-w}\,dz=G_a(w)+\frac{\ln a}{(1-\gamma)(1-w)}.
$$

We now let $\eps\to 0$ in the contour integration. The contribution
from the $\eps$-half circle tends to 0 and we get
$$
\frac{1}{2\pi}\int_{-\pi}^{\pi}\frac{G_a((k+\tfrac12)e^{i\t})}{(k+\tfrac12)e^{i\t}-w}(k+\tfrac12)e^{i\t}\,d\t
-\int_{-k-\tfrac12}^0 \frac{\rho_a(-t)}{t-w}\,dt=G_a(w)+\frac{\ln a}{(1-\gamma)(1-w)}.
 $$

Finally, letting $k\to\infty$ we get \eqref{eq:Garep}, leaving the
details to the reader. Clearly, $\rho_a\ge 0$ if and only if $d_a$
defined in \eqref{eq:darep} is non-negative. It follows that
$G_a(x+1)$ is a Stieltjes function for $a \ge 1$ but not for $0<a<1$,
since in the latter case $\Im G_a(1+iy)>0$ for $y>0$ sufficiently small.
$\quad\square$

\begin{rem} {\rm The integral representation in Theorem~\ref{thm:Ga}
    was established in \cite[(6)]{B:P1} in the case of $a=1$. Since
$$
G_a(z)=G_1(z)+\ln(a)\frac{z}{\log\Gamma(z+1)},
$$
the formula for $G_a$ can be deduced from the formula for $G_1$ and
the following formula
\begin{equation}\label{eq:zlogGamma}
\frac{z}{\log\Gamma(z+1)}=\frac{1}{(1-\gamma)(z-1)}+\int_0^\infty 
\frac{\tau(t)dt}{z+t},\quad z\in\mathcal A\setminus\{1\},
\end{equation}
where
\begin{equation}\label{eq:tau}
\tau(t)=\frac{(k-1)t}{(\ln
  |\Gamma(1-t)|)^2+((k-1)\pi)^2}\quad \mbox{for}
\quad t\in \left]k-1,k\right[,\quad k=1,2,\ldots.
\end{equation}
}
\end{rem}

\end{document}